\documentclass{amsart}
\usepackage{graphicx,url} % Required for inserting images

\usepackage{thmtools} %For renaming theorem 1.3
\usepackage{thm-restate}

\title{Failed Vertex Covering Graph Reconfiguration}

\author{Dennis Hall}
\address{Department of Mathematics, Angelo State University, San Angelo, TX 76904, USA.}
\email{dennis.hall@angelo.edu}

\author{Rogelio Tovar}
\address{Department of Mathematics, Angelo State University, San Angelo, TX 76904, USA.}
\email{rtovar3@angelo.edu}

\date{June 2026}

\usepackage{amsmath}
\usepackage{amssymb}

\usepackage{soul} % For highlighting

\usepackage{comment} %to comment out large blocks

\usepackage{graphicx}
\usepackage{subcaption}
\usepackage[export]{adjustbox}

\usepackage{pgfplots}
\pgfplotsset{compat=1.18}
\usepackage{tikz}
\usetikzlibrary{calc, intersections}
\usepackage{multicol}

\usepackage{color}
\definecolor{darkred}{rgb}{0.7,0,0} % darkred color

\newtheorem{thm}{Theorem}[section]
\newtheorem{theorem}{Theorem}[section]
\newtheorem{proposition}[theorem]{Proposition}
\newtheorem{corollary}[theorem]{Corollary}

\newtheorem{definition}[theorem]{Definition}

\newtheorem{lemma}[theorem]{Lemma}

\theoremstyle{remark}
\newtheorem*{aistatement}{AI Disclosure Statement}

\newcommand{\fzftar}[1]{\mathcal{FZF}(#1)}%^${{\tiny\rm TAR}}
\newcommand{\fvctar}[1]{\mathcal{FVC}(#1)}%^${{\tiny\rm TAR}}

\newcommand{\definitioncolor}{\color{darkred}}
\newcommand{\define}[1]{{\definitioncolor\emph{#1}}} % emphasis of a definition

\definecolor{ASU}{HTML}{4b14f1}
\definecolor{ASUA}{HTML}{942a22}

  \setstcolor{ASU} 
  \setstcolor{ASUA}

\begin{document}

\begin{abstract}
    A failed vertex cover of a graph is a vertex set that leaves at least one edge uncovered. The associated token addition and removal (TAR) graph has the failed vertex covers as vertices, with two sets adjacent when they differ by one vertex. We characterize when two finite graphs have isomorphic failed vertex cover TAR graphs, allowing disconnected graphs, loops, and parallel edges.
\end{abstract}

\maketitle

\section{Introduction}{
As discussed by Curtis, Flagg, and Hogben in \cite{XTAR-YTAR}, graph parameters such as vertex covering number, zero forcing number, and power domination number can be abstracted with the concept of $X$-set parameters. For the parameters considered here, the families of $X$-sets are closed under taking supersets, while the corresponding failed-set families are closed under taking subsets. This is the basis for the analogous concept, $Y$-sets. Token addition and removal (TAR) graphs are a tool used to map out these types of parameters for a given graph. The following isomorphism theorem for $Y$-sets was established, where $\mathfrak{Y}(G)$ is the $Y$-set TAR graph of $G$.

\begin{thm}[Isomorphism Theorem for robust $Y$-set parameters]
\cite[Theorem 8.13]{XTAR-YTAR} 
\label{thm:IsomorphismRobust}
Let $Y$ be a robust $Y$-set parameter and let $G$ and $G'$ be base graphs such that $\mathfrak{Y}(G)\cong\mathfrak{Y}(G')$. If $Y(K_1)=1$ or $G$ and $G'$ have no isolated vertices, then $G$ and $G'$ have the same order and there is a relabeling of the vertices of $G'$ such that $G$ and $G'$ have exactly the same $Y$-sets.
\end{thm}

Note that this theorem requires the {\it robustness} hypothesis. Many common $Y$-set parameters, such as failed zero forcing and failed vertex covering, are not robust, though a similar theorem can be found if we require that the graphs are connected. In \cite{ZFpaper}, an analogous theorem is established for failed zero forcing.

\begin{thm}[Isomorphism Theorem for Failed Zero Forcing]
\cite[Theorem 4.17]{ZFpaper}
\label{maintheorem}
Let $G$ and $G'$ be connected graphs on at least two vertices such that $\fzftar{G} \cong \fzftar{G'}$. Then $G$ and $G'$ have the same order and there is a relabeling of the vertices of $G'$ so that $G$ and $G'$ have exactly the same failed zero forcing sets.
\end{thm}

For failed zero forcing, the TAR graph loses much of the information. For failed vertex covers, however, most of the information can be recovered. In the following theorem, we characterize those graphs that have isomorphic FVC-TAR graphs. The reduced graph $R(G)$ is obtained by replacing each parallel-edge class with a single edge and deleting each nonloop edge incident to a looped vertex. The graph $L_l$ is a graph consisting of $l$ looped vertices and no other edges. The graph $iK_1$ is the graph consisting of $i$ isolated vertices and no edges. All of these terms are defined more precisely in Section 5.

\begin{restatable}[Isomorphism Theorem for Failed Vertex Covering in Graphs]{thm}{maintheorem}
    Let $G$ and $G'$ be graphs. Then $\fvctar{G}\cong\fvctar{G'}$ if and only if one of the following holds.
    \begin{enumerate}
        \item Both $G$ and $G'$ are edgeless.
        \item Both $G$ and $G'$ have at least one edge and $R(G)\cong R(G')$.
        \item After possibly interchanging $G$ and $G'$, there exist integers $l\geq 1$ and $i\geq 0$ such that $R(G)\cong K_{1,l}\oplus iK_1$ and $R(G')\cong L_l\oplus iK_1$.
    \end{enumerate}
    In particular, if $G$ and $G'$ are both connected and simple, then $G\cong G'$ if and only if $\fvctar{G}\cong\fvctar{G'}$.
\end{restatable}

Section 2 introduces the necessary terminology and preliminary results, Section 3 develops the theorem for connected simple graphs, Section 4 studies disconnected graphs, and Section 5 examines graphs with loops and parallel edges and proves the main theorem.

\subsection*{Acknowledgments}
Research was sponsored by the Air Force Office of Scientific Research and was accomplished under Grant Number W911NF-25-1-0110. The views and conclusions contained in this document are those of the authors and should not be interpreted as representing the official policies, either expressed or implied, of the Army Research Office or the U.S. Government. The U.S. Government is authorized to reproduce and distribute reprints for Government purposes notwithstanding any copyright notation herein.

\section{Preliminaries}

Here, we establish notation for some terms in graph and reconfiguration theory. Throughout the paper, we follow \cite{Graph Theory} for graph-theoretic terminology and \cite{XTAR-YTAR} for reconfiguration terminology.
Though we provide an overview of the most-used definitions below, the reader is referred to these citations for a more thorough treatment. All graphs in Sections 2--4 are simple.

\subsection{Graph Theory Terminology} A \define{graph} $G$ is a collection of vertices $V(G)$ with a set of edges $E(G)$, which are unordered pairs of vertices. All graphs in this paper are finite, and all base graphs have a nonempty vertex set. The number of vertices adjacent to a given vertex $v$ is the \define{degree} of $v$, denoted $\deg(v)$. The maximum degree among all the vertices in a graph $G$ is denoted by $\Delta(G)$ and the minimum degree is $\delta(G)$.

For sets $A$ and $B$, we write $A\ominus B=(A\setminus B)\cup (B\setminus A)$ for their symmetric difference. The graph $G\setminus L$ is obtained by deleting the vertices in $L$ and all incident edges. $G\sqcup H$ denotes disjoint union.

\begin{definition}[Vertex Cover]
    Given a graph $G$, a \define{vertex cover} is a set of vertices $C\subseteq V(G)$ such that every edge in $E(G)$ is incident to at least one vertex in $C$. The minimum cardinality of such a vertex cover is the \define{vertex covering number} of $G$, denoted $\tau(G)$. A set of vertices that is not a vertex cover is a \define{failed vertex cover}.
\end{definition}

Observe that one vertex cover for the graph $C_5$, as shown in Figure \ref{C_5 Cover}, is $\{1,3,4\}$. Note that the vertex set $\{2,5\}$ is a {\it failed} vertex cover since neither $2$ nor $5$ is incident to the edge $(3,4)$. Indeed, there is no vertex cover of size $2$ in $C_5$ and thus $\tau(C_5)=3$.

\begin{figure}[h!]
    \begin{tikzpicture}[
    scale=0.5,
  vertex/.style={circle,draw=black,fill=white,minimum size=8mm,inner sep=0pt,font=\small},
  selected/.style={vertex,fill=black!25},
  edge/.style={black,line width=0.6pt}]
  \node[selected] (v1) at (0,2.7) {$1$};
  \node[vertex]   (v2) at (2.57,0.83) {$2$};
  \node[selected] (v3) at (1.59,-2.2) {$3$};
  \node[selected] (v4) at (-1.59,-2.2) {$4$};
  \node[vertex]   (v5) at (-2.57,0.83) {$5$};
  \foreach \a/\b in {1/2,2/3,3/4,4/5,5/1} \draw[edge] (v\a)--(v\b);
\end{tikzpicture}
\caption{The cycle $C_5$ with a vertex cover shaded.}\label{C_5  Cover}
\end{figure}
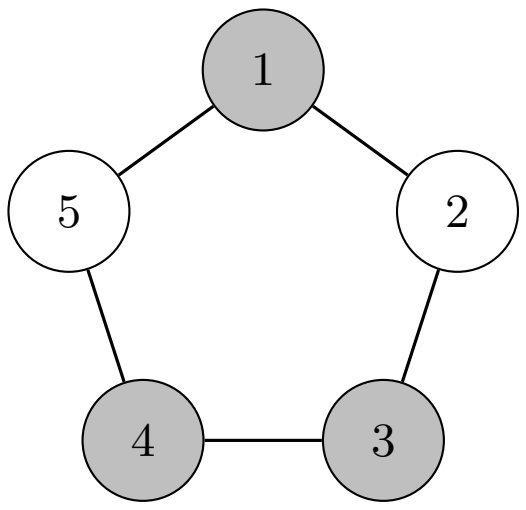

\begin{definition} A \define{star graph}, denoted $S_n$ for $n\in \mathbb{Z}^+$, is a simple connected graph with a central vertex $v$ such that every other vertex is adjacent only to $v$. Equivalently, $S_n=K_{1,n}$ has $n+1$ vertices and $n$ leaves.
\end{definition}

A well-known result (see, e.g., \cite{starcover}) is that a star is the only connected graph with $\tau(G)=1$.

\begin{proposition}\label{starunique}
    A connected graph $G$ has $\tau(G)=1$ if and only if $G$ is a star.
\end{proposition}

\begin{definition}The \define{failed vertex cover token addition and removal (FVC-TAR) graph} of a graph $G=(V,E)$, denoted by $\fvctar{G}$, is the graph with vertices the subsets of $V$ that are failed vertex covers on $G$ with edges so that two vertices are adjacent if and only if they differ by exactly one element. This graph is an induced subgraph of a hypercube.
\end{definition}

In Figure \ref{F(C_4) Labeled}, we show $\fvctar{C_5}$. Observe that every vertex set of size less than three is a failed vertex cover of $C_5$, and only certain size-3 sets are failed vertex covers. Note that the isomorphism in the hypothesis of Lemma \ref{thm:main} is between unlabeled copies of the FVC-TAR graphs. Though recovering the labels themselves is often impossible, we are able to find a valid relabeling. The relationship between TAR graph isomorphisms and vertex relabelings is explored more thoroughly in \cite{XTAR-YTAR}.

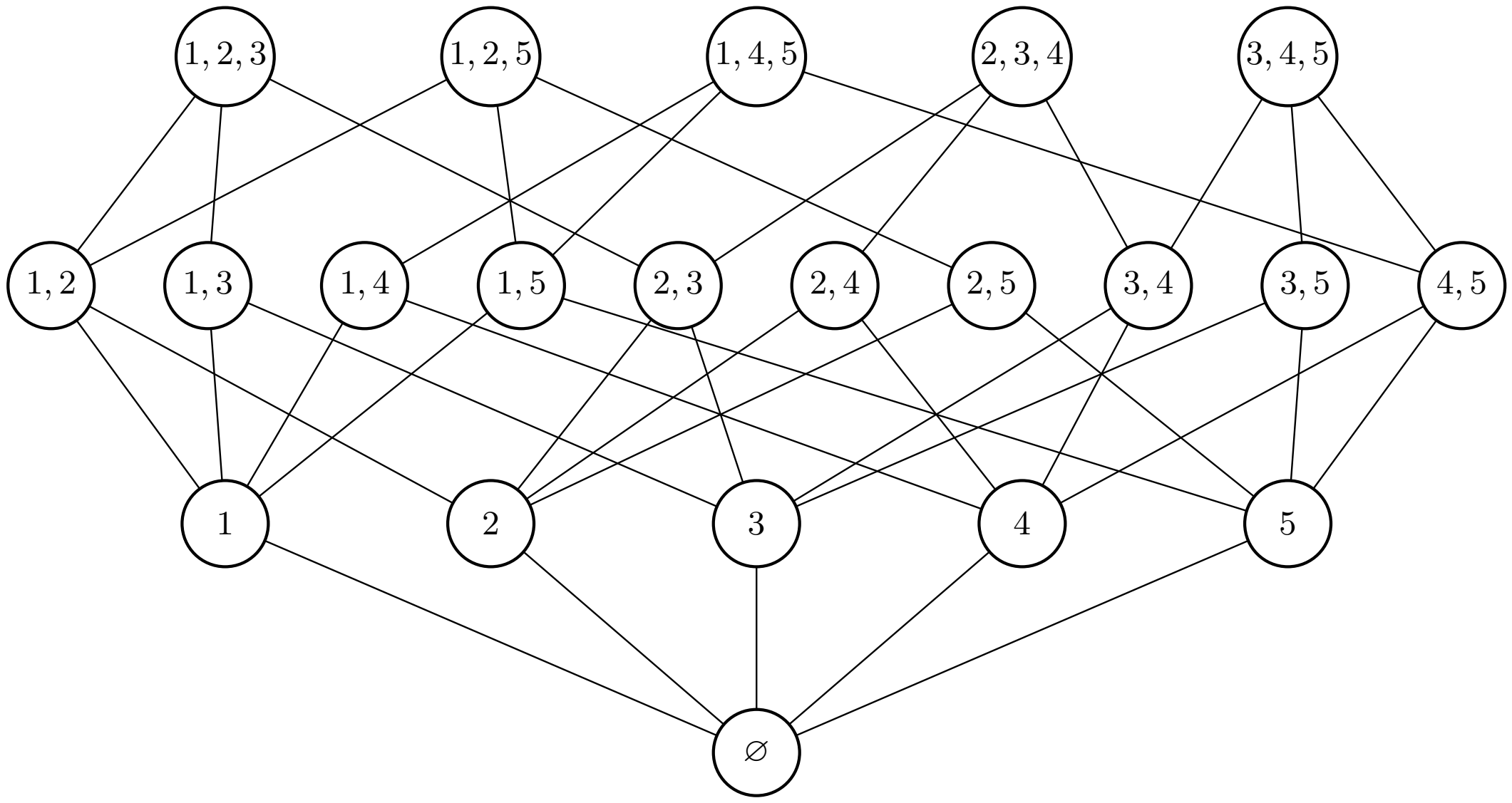
\begin{figure}[h!]
    \begin{tikzpicture}[
  x=0.85cm,y=0.85cm,
  vertex/.style={circle,draw=black,fill=white,line width=0.8pt,minimum size=8mm,inner sep=1pt,font=\small},
  edge/.style={draw=black,line width=0.45pt}
]
% C5 has edges 12,23,34,45,51.
% Vertices are all subsets that fail to cover at least one of these edges.

% Coordinates for edge endpoints
\coordinate (v123) at (-5.8000,8.1000);
\coordinate (v125) at (-2.9000,8.1000);
\coordinate (v145) at (0.0000,8.1000);
\coordinate (v234) at (2.9000,8.1000);
\coordinate (v345) at (5.8000,8.1000);
\coordinate (v12) at (-7.7000,5.6000);
\coordinate (v13) at (-5.9889,5.6000);
\coordinate (v14) at (-4.2778,5.6000);
\coordinate (v15) at (-2.5667,5.6000);
\coordinate (v23) at (-0.8556,5.6000);
\coordinate (v24) at (0.8556,5.6000);
\coordinate (v25) at (2.5667,5.6000);
\coordinate (v34) at (4.2778,5.6000);
\coordinate (v35) at (5.9889,5.6000);
\coordinate (v45) at (7.7000,5.6000);
\coordinate (v1) at (-5.8000,3.0000);
\coordinate (v2) at (-2.9000,3.0000);
\coordinate (v3) at (0.0000,3.0000);
\coordinate (v4) at (2.9000,3.0000);
\coordinate (v5) at (5.8000,3.0000);
\coordinate (empty) at (0.0000,0.5000);

% Draw edges behind vertices
\draw[edge] (empty) -- (v1);
\draw[edge] (empty) -- (v2);
\draw[edge] (empty) -- (v3);
\draw[edge] (empty) -- (v4);
\draw[edge] (empty) -- (v5);
\draw[edge] (v1) -- (v12);
\draw[edge] (v1) -- (v13);
\draw[edge] (v1) -- (v14);
\draw[edge] (v1) -- (v15);
\draw[edge] (v2) -- (v12);
\draw[edge] (v2) -- (v23);
\draw[edge] (v2) -- (v24);
\draw[edge] (v2) -- (v25);
\draw[edge] (v3) -- (v13);
\draw[edge] (v3) -- (v23);
\draw[edge] (v3) -- (v34);
\draw[edge] (v3) -- (v35);
\draw[edge] (v4) -- (v14);
\draw[edge] (v4) -- (v24);
\draw[edge] (v4) -- (v34);
\draw[edge] (v4) -- (v45);
\draw[edge] (v5) -- (v15);
\draw[edge] (v5) -- (v25);
\draw[edge] (v5) -- (v35);
\draw[edge] (v5) -- (v45);
\draw[edge] (v12) -- (v123);
\draw[edge] (v12) -- (v125);
\draw[edge] (v13) -- (v123);
\draw[edge] (v14) -- (v145);
\draw[edge] (v15) -- (v125);
\draw[edge] (v15) -- (v145);
\draw[edge] (v23) -- (v123);
\draw[edge] (v23) -- (v234);
\draw[edge] (v24) -- (v234);
\draw[edge] (v25) -- (v125);
\draw[edge] (v34) -- (v234);
\draw[edge] (v34) -- (v345);
\draw[edge] (v35) -- (v345);
\draw[edge] (v45) -- (v145);
\draw[edge] (v45) -- (v345);

% Level 3
\node[vertex] (v123) at (-5.8000,8.1000) {$1,2,3$};
\node[vertex] (v125) at (-2.9000,8.1000) {$1,2,5$};
\node[vertex] (v145) at (0.0000,8.1000) {$1,4,5$};
\node[vertex] (v234) at (2.9000,8.1000) {$2,3,4$};
\node[vertex] (v345) at (5.8000,8.1000) {$3,4,5$};

% Level 2
\node[vertex] (v12) at (-7.7000,5.6000) {$1,2$};
\node[vertex] (v13) at (-5.9889,5.6000) {$1,3$};
\node[vertex] (v14) at (-4.2778,5.6000) {$1,4$};
\node[vertex] (v15) at (-2.5667,5.6000) {$1,5$};
\node[vertex] (v23) at (-0.8556,5.6000) {$2,3$};
\node[vertex] (v24) at (0.8556,5.6000) {$2,4$};
\node[vertex] (v25) at (2.5667,5.6000) {$2,5$};
\node[vertex] (v34) at (4.2778,5.6000) {$3,4$};
\node[vertex] (v35) at (5.9889,5.6000) {$3,5$};
\node[vertex] (v45) at (7.7000,5.6000) {$4,5$};

% Level 1
\node[vertex] (v1) at (-5.8000,3.0000) {$1$};
\node[vertex] (v2) at (-2.9000,3.0000) {$2$};
\node[vertex] (v3) at (0.0000,3.0000) {$3$};
\node[vertex] (v4) at (2.9000,3.0000) {$4$};
\node[vertex] (v5) at (5.8000,3.0000) {$5$};

% Level 0
\node[vertex] (empty) at (0.0000,0.5000) {$\varnothing$};

% TAR edges: one element added or removed
\end{tikzpicture}
\caption{The labeled FVC-TAR graph of $C_5$.}\label{F(C_4) Labeled}
\end{figure}
}

\section{Results on Simple Connected Graphs}

In this section, we provide results on simple connected graphs. All graphs in this section are simple.

\begin{proposition}\label{star unique}
    If $\fvctar{S_m}\cong\fvctar{S_n}$ then $m=n$.
\end{proposition}
\begin{proof}
Let $S_k$ be a star graph with central vertex $v$. The vertex covers of $S_k$ are exactly those $2^k$ sets that contain $v$ and the two sets that contain all of $V(S_k)\setminus\{v\}$. Since there is one set that contains both $v$ and $V(S_k)\setminus\{v\}$, there are a total of $2^k+2-1=2^k+1$ vertex covers of $S_k$ and thus $2^{k+1}-(2^k+1)=2^k-1$ failed vertex covers. Then $|V(\fvctar{S_k})|=2^k-1$. Applying this count to the assumed isomorphism, since $\fvctar{S_n}\cong \fvctar{S_m}$, it follows that \[2^m-1=|V(\fvctar{S_m})|=|V(\fvctar{S_n})|=2^n-1\] and therefore $m=n$.

\end{proof}

\begin{lemma} \label{vertex count}
    Let $G$ be a connected graph. If $\fvctar{G}\cong\fvctar{S_n}$ then either $G\cong S_n$ or $|V(G)|=n$. 
\end{lemma}
\begin{proof}
    We assume $G\not\cong S_n$ and show that $|V(G)|=n$. By Proposition \ref{star unique}, this implies that $G\not\cong S_m$ for any $m$. We now consider the case in which $G$ is not a star. Since $G$ is connected, has no isolated vertices, and is not a star, it has at least three vertices. We know that 
    $\tau(G) \geq 2$ by Proposition \ref{starunique}.  This means that no singleton is a vertex cover of $G$. Thus, every singleton is represented by a vertex in $\fvctar{G}$. The vertex of $\fvctar{G}$ corresponding to the empty set is adjacent to every singleton and therefore has degree $|V(G)|$. Because $\fvctar{G}$ toggles the vertices of $G$, and at most $|V(G)|$ elements can be toggled, it follows that $\Delta(\fvctar{G}) = |V(G)|$. 
    
    %\Dennis{If $n=1$, then $S_1=K_2$ and $\fvctar{S_1}$ consists only of the vertex $\varnothing$. Its maximum degree is therefore zero, contradicting the previous paragraph which established that $\Delta(\mathcal{F}(G)) = |V(G)|\geq 2$. Thus $n\geq 2$. Or something like this. $n=1$ doesn't work.}

    Let us now consider $\fvctar{S_n}$. Note that if $n=1$, then $\fvctar{S_1}=K_1$, contradicting $\Delta(\fvctar{G})= |V(G)|\geq 3$. Thus $n\geq 2$. We know that because the vertex covering number of $S_n$ is 1, there exists a vertex in $S_n$, $c$, such that any subset that contains $c$ will not be in $V(\fvctar{S_n})$. Thus, for all $v \in V(\fvctar{S_n}),$ it holds that $ \deg(v) \leq n$, which is the same as
    $\Delta(\fvctar{S_n})\leq n$. We observe that the vertex in $\fvctar{S_n}$ that represents the empty set will have edges connecting to all of the vertices in $\fvctar{S_n}$ corresponding to the singletons, of which there are $n$. Thus, the vertex in $\fvctar{S_n}$ corresponding to the empty set has degree $n$.

    Because $\fvctar{G} \cong \fvctar{S_n}$, it follows that $\Delta(\fvctar{G})=\Delta(\fvctar{S_n}) \implies |V(G)|=n$, and so the conclusion follows.
\end{proof}

\begin{lemma} \label{subgraph lemma}
    If $G$ is a subgraph of a graph $H$, then $\fvctar{G}$ is an induced subgraph of $\fvctar{H}$.
\end{lemma}
\begin{proof}
    Let $A$ be a failed vertex cover of $G$. Since $A$ leaves some edge of $G$ uncovered and that edge also belongs to $H$, the set $A$ is a failed vertex cover of $H$. Consider the subgraph of $\fvctar{H}$ induced by the collection of failed vertex covers of $G$. Two such sets are adjacent precisely when they differ by one element, which is also their adjacency condition in $\fvctar{G}$. Therefore this induced subgraph is $\fvctar{G}$.
\end{proof}

The following corollary follows immediately.

\begin{corollary}\label{kncorr}
    If $G$ is a graph on $n$ vertices, then $\fvctar{G}$ is an induced subgraph of $\fvctar{K_n}$.
\end{corollary}

\begin{proposition}\label{propstar}
    Let $G$ be a connected graph. If $\fvctar{G}\cong\fvctar{S_n}$ then $G\cong S_n$.
\end{proposition}
\begin{proof}
    Suppose, for the sake of contradiction, that there exists a graph $G$ such that $\fvctar{G} \cong \fvctar{S_n}$, but $G \ncong S_n$. Then $|V(G)|=n$ by Lemma \ref{vertex count}. Because 
    $\fvctar{G} \cong \fvctar{S_n}$, it follows that $|V(\fvctar{G})|=|V(\fvctar{S_n})|$. 
    By Corollary \ref{kncorr}, we know that $\fvctar{G}$ is a subgraph of $\fvctar{K_n}$, and so $|V(\fvctar{K_n})|\geq |V(\fvctar{G})| = |V(\fvctar{S_n})|$. Since the failed vertex covers of $K_n$ are exactly the $2^n-n-1$ vertex sets of cardinality at most $n-2$ and  those of $S_n$ are the $2^n-1$ proper subsets of the $n$ non-central vertices, we have
    \[2^n - n - 1 = |V(\fvctar{K_n})| \geq |V(\fvctar{S_n})| = 2^n-1,\]
    a contradiction for any $n\geq 1$. Therefore, $G\cong S_n$, as desired.
\end{proof}

\begin{proposition}
    For $n\geq 1$, $\fvctar{K_n}$ is an induced subgraph of $\fvctar{S_n}$.
\end{proposition}
\begin{proof}
    The failed vertex covers of $K_n$ are exactly its vertex subsets of size at most $n-2$. Identify $V(K_n)$ with the set of $n$ leaves of $S_n$. Under this identification, every failed vertex cover of $K_n$ is a proper subset of the leaves and hence a failed vertex cover of $S_n$. Adjacency in both TAR graphs is determined by toggling elements. Thus $\fvctar{K_n}$ is an induced subgraph of $\fvctar{S_n}$.
\end{proof}

\begin{corollary}
    If $G$ is a graph on $n\geq 1$ vertices, then $\fvctar{G}$ is an induced subgraph of $\fvctar{S_n}$.
\end{corollary}

    The family of failed vertex covers is closed under taking subsets. Furthermore, if $G$ is a connected graph of order at least two that is not a star, then Proposition \ref{starunique} implies that every singleton is a failed vertex cover. These observations allow us to construct a robust parameter that agrees with failed vertex covering on connected non-star graphs. We now prove a relabeling result for connected simple graphs.

\begin{lemma}
\label{thm:main}
Let $G$ and $G'$ be connected graphs on at least two vertices such that
$\fvctar{G} \cong \fvctar{G'}$. Then $G$ and $G'$ have the same order and
there is a relabeling of the vertices of $G'$ so that $G$ and $G'$ have
exactly the same failed vertex covers.
\end{lemma}
\begin{proof}
    Suppose first that $G\cong S_n$. Since $\fvctar{G'}\cong \fvctar{S_n}$, Proposition \ref{propstar} implies that $G'\cong S_n$, and the conclusion follows.

    Now suppose that $G$ is not a star. By Proposition \ref{propstar}, $G'$ is not a star either. We are not able to apply Theorem \ref{thm:IsomorphismRobust} directly since the failed vertex cover families do not define a robust $Y$-set parameter. We construct an auxiliary parameter $Y^*$ as follows. On a connected non-star graph of order at least two, the $Y^*$-sets are exactly the failed vertex covers. On a star or $K_1$, every vertex subset is a $Y^*$-set. For a graph $H$ with connected components $H_1,\ldots,H_k$, a set $A\subseteq V(H)$ is a $Y^*$-set if and only if $A\cap V(H_i)$ is a $Y^*$-set of $H_i$ for every $i$. Define $Y^*(H)$ to be the maximum cardinality of a $Y^*$-set. 
    
    This construction is invariant under graph isomorphism, and its admissible families are nonempty and closed under taking subsets. Every singleton is a $Y^*$-set on a connected graph of order at least two, which follows from Proposition \ref{starunique} for non-stars and from the definition for stars. Component consistency holds by construction. Thus $Y^*$ is robust with $Y^*(K_1)=1$. Its TAR graphs on $G$ and $G'$ are precisely $\fvctar{G}$ and $\fvctar{G'}$. The hypotheses of Theorem \ref{thm:IsomorphismRobust} are therefore satisfied, so the theorem gives the desired equality of orders and relabeling.
\end{proof}

Curtis, Flagg, and Hogben \cite{XTAR-YTAR} proved that a simple graph is determined up to isomorphism by its vertex cover TAR graph. The following corollary establishes the analogous result for failed vertex cover TAR graphs of connected simple graphs.

\begin{corollary}\label{connectediso}
    Let $G$ and $G'$ be connected simple graphs. Then $\fvctar{G}\cong\fvctar{G'}$ if and only if $G\cong G'$.
\end{corollary}
\begin{proof}
    The reverse direction is immediate, and therefore we assume $\fvctar{G}\cong\fvctar{G'}$. If $G\cong K_1$, then $\fvctar{G}$ has no vertices and thus $\fvctar{G'}$ also has no vertices. Since the empty set would be a failed vertex cover if $G'$ had an edge, it follows that $G'$ is edgeless. Connectedness therefore implies that $G'\cong K_1\cong G$, as desired. We thus assume that both graphs have order at least two.

    By Lemma \ref{thm:main}, we may relabel the vertices of $G'$ so that $G$ and $G'$ have the same vertex set $V$ and exactly the same failed vertex covers. Let $u,v\in V$ be distinct. Since $G$ is simple, the only possible edge left uncovered by $V\setminus\{u,v\}$ is $uv$. Therefore $uv\in E(G)$ if and only if $V\setminus\{u,v\}$ is a failed vertex cover of $G$. The same equivalence holds for $G'$. Because the two graphs have exactly the same failed vertex covers, it follows that $uv\in E(G)$ if and only if $uv\in E(G')$. Thus $E(G)=E(G')$ under this relabeling and hence $G\cong G'$.
\end{proof}

\section{Disconnected Graphs}
  Following an observation made in \cite{ZFpaper}, we note that the graph $\fvctar{G}$ can be viewed as a subgraph of the $n$-dimensional hypercube $Q_n$, where $n=|V(G)|$, the vertices correspond to subsets of $V(G)$, and edges correspond to symmetric differences of size one. Further, we have the following lemma, the downward-closed analogue of the distance property in \cite[Proposition 2.5]{powerdom}. The \define{distance} $\textrm{d}_G(u,v)$ between two vertices $u$ and $v$ of a finite connected graph $G$ is the minimum length of the paths connecting them. In this section, all graphs are simple.

\begin{lemma}\label{symmdistance} \cite[Proposition 2.5]{powerdom} Let $G$ have at least one edge. Then $\fvctar{G}$ is connected and \[\textrm{d}_{\fvctar{G}}(S,T)=|S\ominus T|\quad \textrm{ for }\quad S,T\in V(\fvctar{G}).\] In particular, $\textrm{d}_{\fvctar{G}}(\varnothing,S)=|S|$.
\end{lemma}

We wish to show that we are able to identify the connected components of $G$ within its unlabeled FVC-TAR graph without relying on the original graph. To begin, we must first find which element corresponds to the empty set. The \define{eccentricity} $\epsilon(v)$ of a vertex $v$ in a connected graph $G$ is the maximum distance between $v$ and any other vertex $u$ of $G$. The \define{radius} $r(G)$ is the minimum eccentricity found across all vertices in the graph. The \define{center} is the set of vertices with eccentricity equal to $r(G)$.

\begin{lemma}\label{findcenter}
    Let $G$ be a simple graph with no isolated vertices and at least one edge. Then the empty set is the unique central vertex of $\fvctar{G}$ and $|V(G)|=r(\fvctar{G})+2$.
\end{lemma}
\begin{proof}
    Set $m=|V(G)|$. A failed vertex cover must omit both endpoints of at least one edge and thus has at most $m-2$ elements. For every edge $e$, the set $V(G)\setminus e$ is a failed vertex cover that achieves this maximum size. Observing that the symmetric difference between the empty set and these sets is also maximum, it follows from Lemma \ref{symmdistance} that $\epsilon_{\fvctar{G}}(\varnothing)=m-2$.

    We now show that every other vertex in $\fvctar{G}$ has eccentricity greater than that of the empty set. Let $S\neq\varnothing$ be a vertex of $\fvctar{G}$. If $S$ contains an edge, its complement $T=V(G)\setminus S$ misses that edge and is thus represented by a vertex of $\fvctar{G}$. Since $\textrm{d}_{\fvctar{G}}(S,T)=|S\ominus T|=m$ and no symmetric difference can have size greater than $m$, it follows that $\epsilon_{\fvctar{G}}(S)=m$.
    
    On the other hand, if $S$ does not contain an edge, choose $v\in S$. Since $G$ has no isolated vertices, it follows that $v$ has a neighbor $w\not\in S$. Then $T=V(G)\setminus (S\cup \{w\})$ misses the edge $vw$ and is thus a failed vertex cover. Then $\textrm{d}_{\fvctar{G}}(S,T)=|S\ominus T|=|V(G)-\{w\}|=m-1$. The only set that can have greater distance from $S$ would be $S^C=V(G)\setminus S$. However, since $S$ contains no edges, it follows that $S^C$ intersects every edge and is thus a vertex cover. Then $\epsilon_{\fvctar{G}}(S)=m-1$.

    Therefore every nonempty failed vertex cover has eccentricity strictly larger than that of $\varnothing$ and every other failed vertex cover has eccentricity at least $m-1$. Hence the center contains only $\varnothing$ and $r(\fvctar{G})=m-2$, that is, $|V(G)|=r(\fvctar{G})+2$.
\end{proof}

\begin{proposition} \label{disconnected main theorem}
    Let $G$ and $G'$ be simple graphs with no isolated vertices. Then $\fvctar{G}\cong\fvctar{G'}$ if and only if $G\cong G'$. 
\end{proposition}
\begin{proof}
    The reverse implication is immediate. Since $G$ and $G'$ have no isolated vertices, they must have at least one edge. Then, by Lemma \ref{findcenter}, $|V(G)| = r(\fvctar{G})+2 = r(\fvctar{G'})+2 = |V(G')|$. If $G$ is a star with central vertex $c$, then $\{c\}$ is a vertex cover, and the degree of the empty-set vertex in $\fvctar{G}$ is less than $m=|V(G)|$. The same holds, then, for $G'$, and we can find $c'\in V(G')$ so that $\{c'\}$ is a vertex cover in $G'$. Since $G'$ has no isolated vertices, every other vertex is adjacent to $c'$, and every edge is incident to $c'$. Thus $G'$ is also a star. Now, Proposition \ref{star unique} applies and so $G\cong G'$.

     We now assume that neither $G$ nor $G'$ is a star. Every singleton is thus a failed vertex cover in both graphs. Let $\phi: \fvctar{G}\to\fvctar{G'}$ be an isomorphism and define $f: V(G)\to V(G')$ by $f(v) = x$ when $\phi(\{v\})=\{x\}$. Since every singleton is a failed vertex cover, the map $f$ is a bijection. For $S\subseteq V(G)$, define $f(S) = \{f(v) : v\in S\}$.
     
     Now, let $v\in S$. Then, by Lemma \ref{symmdistance}, $d_{\fvctar{G}}(\{v\},S)=d_{\fvctar{G}}(\varnothing,S)-1$. As $\phi$ is an isomorphism between $\fvctar{G}$ and $\fvctar{G'}$, and $\phi(\varnothing)=\varnothing$, it follows that $d_{\fvctar{G'}}(\phi(\{v\}),\phi(S))=d_{\fvctar{G'}}(\varnothing,\phi(S))-1$. By the definition of $f$, then, $d_{\fvctar{G'}}(\{f(v)\},\phi(S))=d_{\fvctar{G'}}(\varnothing,\phi(S))-1$. Applying the symmetric difference equality from Lemma \ref{symmdistance}, this implies that $|\{f(v)\} \ominus \phi(S)|=|\varnothing\ominus\phi(S)|-1=|\phi(S)|-1$. This occurs exactly when $f(v)\in\phi(S)$. The converse at each step also holds and thus $v\in S$ if and only if $f(v)\in \phi(S)$. Therefore, since $f$ is surjective, we have $\phi(S)=f(S)$.

     Finally, let $u,v$ be distinct vertices in $V(G)$. Then $uv\in E(G)$ exactly when $V(G)\setminus\{u,v\}$ is a failed vertex cover. Since $f$ is a bijection, $f(V(G)\setminus\{u,v\})=V(G')\setminus\{f(u),f(v)\}$. Because $\phi(S)=f(S)$, it follows that $\phi(V(G)\setminus\{u,v\})=V(G')\setminus\{f(u),f(v)\}$. As $\phi$ maps failed vertex covers bijectively to failed vertex covers, we have that $V(G')\setminus\{f(u),f(v)\}$ is a failed vertex cover, which means $f(u)f(v)\in E(G')$. Therefore $f$ is a graph isomorphism.
\end{proof}

\section{Graphs with Loops and Parallel Edges}

Now, we will extend this idea to graphs that may have parallel edges and loops. In this section, graphs need not be connected, and isolated vertices are allowed. When deleting vertices leaves the empty graph, we take its vertex covering number to be zero.

For $u,v\in V$, let $\mu_G(u,v)$ be the edge multiplicity between $u$ and $v$. Then $\mu_G$ is symmetric and, if a vertex $u$ has no loops, then $\mu_G(u,u)=0$. The set of all edges between $u$ and $v$ is a \define{parallel-edge class}. If $\mu_G(u,v)=1$, then the set containing exactly the edge between $u$ and $v$ is a \define{trivial parallel-edge class} or \define{single edge}. In particular, multiple loops at one vertex form a single class.

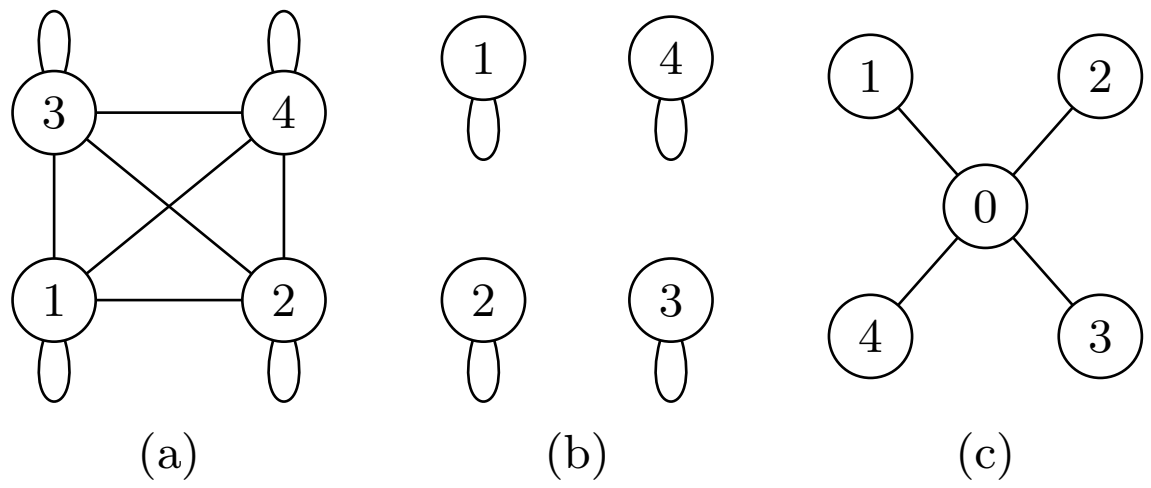
\begin{figure}[h!tbp]
  \centering
  \begin{tikzpicture}[
    line width=0.5pt,
    vertex/.style={
      circle, draw, fill=white, minimum size=5.5mm,
      inner sep=0pt, outer sep=0pt, font=\small
    },
    every loop/.style={min distance=5.6mm}
  ]
    % (a) Second graph in the top row.
    \begin{scope}
      \node[vertex] (a1) at (-0.76,-0.62) {$1$};
      \node[vertex] (a2) at ( 0.76,-0.62) {$2$};
      \node[vertex] (a3) at (-0.76, 0.62) {$3$};
      \node[vertex] (a4) at ( 0.76, 0.62) {$4$};
      \foreach \source/\target in {1/2, 1/3, 1/4, 2/3, 2/4, 3/4}
        \draw (a\source) -- (a\target);
      \path (a3) edge[loop above] (a3)
            (a4) edge[loop above] (a4)
            (a1) edge[loop below] (a1)
            (a2) edge[loop below] (a2);
      \node[font=\small, inner sep=0pt] at (0,-1.65) {(a)};
    \end{scope}

    % (b) Third graph in the top row.
    \begin{scope}[xshift=2.7cm]
      \node[vertex] (b1) at (-0.62, 0.98) {$1$};
      \node[vertex] (b4) at ( 0.62, 0.98) {$4$};
      \node[vertex] (b2) at (-0.62,-0.62) {$2$};
      \node[vertex] (b3) at ( 0.62,-0.62) {$3$};
      \foreach \index in {1,2,3,4}
        \path (b\index) edge[loop below] (b\index);
      \node[font=\small, inner sep=0pt] at (0,-1.65) {(b)};
    \end{scope}

    % (c) Graph in the bottom row.
    \begin{scope}[xshift=5.4cm]
      \node[vertex] (c0) at (0,0) {$0$};
      \node[vertex] (c1) at (-0.76, 0.86) {$1$};
      \node[vertex] (c2) at ( 0.76, 0.86) {$2$};
      \node[vertex] (c3) at ( 0.76,-0.86) {$3$};
      \node[vertex] (c4) at (-0.76,-0.86) {$4$};
      \foreach \index in {1,2,3,4}
        \draw (c0) -- (c\index);
      \node[font=\small, inner sep=0pt] at (0,-1.65) {(c)};
    \end{scope}
  \end{tikzpicture}
  \caption{(a)~$K_4$ with a loop at each vertex.
    (b)~Four disjoint loops. (c)~The star $K_{1,4}$.}
  \label{fig:three-small-graphs}
\end{figure}

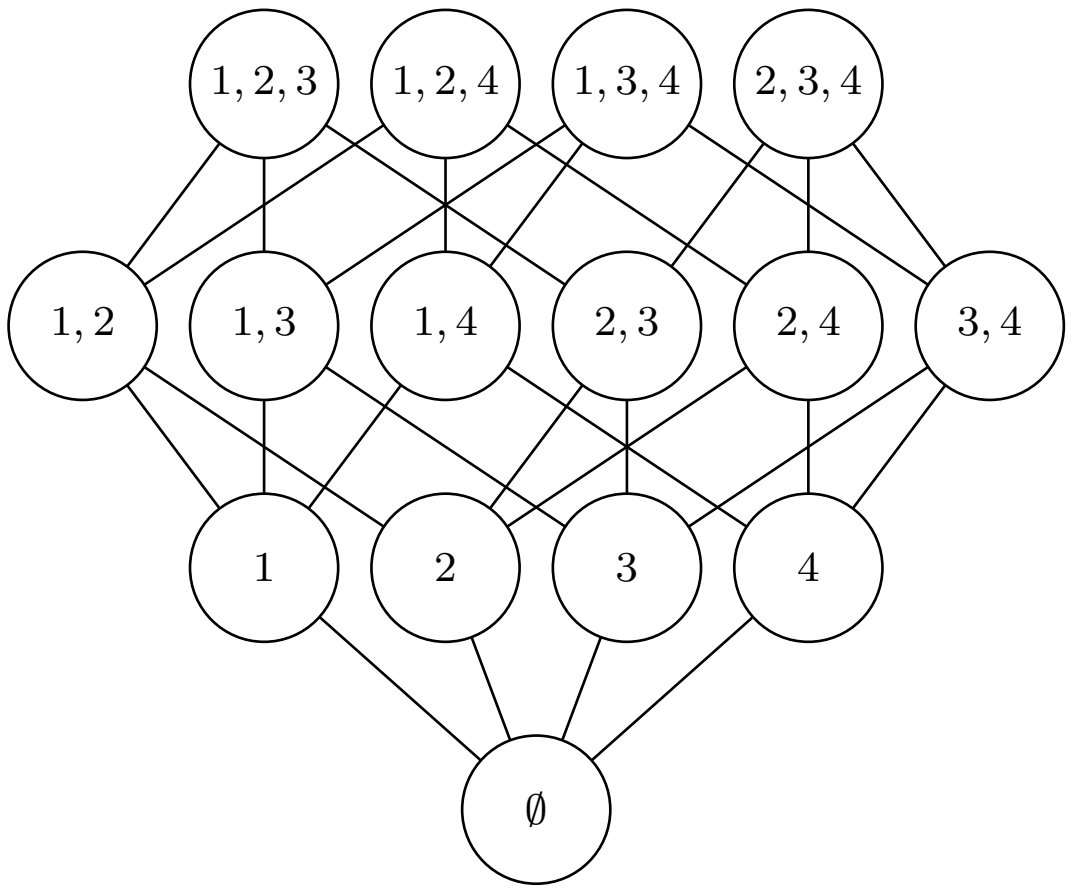
\begin{figure}[tbp]
  \centering
  \begin{tikzpicture}[
    x=12mm, y=16mm, line width=0.5pt,
    subset/.style={
      circle, draw, fill=white, minimum size=9.8mm,
      inner sep=0.5pt, outer sep=0pt, font=\footnotesize
    }
  ]
    \node[subset] (s0) at (0,0) {$\emptyset$};
    \foreach \index/\column in {1/-1.5, 2/-0.5, 3/0.5, 4/1.5}
      \node[subset] (s\index) at (\column,1) {$\index$};
    \foreach \id/\column/\elements in {
      12/-2.5/{1,2}, 13/-1.5/{1,3}, 14/-0.5/{1,4},
      23/0.5/{2,3}, 24/1.5/{2,4}, 34/2.5/{3,4}}
      \node[subset] (s\id) at (\column,2) {$\elements$};
    \foreach \id/\column/\elements in {
      123/-1.5/{1,2,3}, 124/-0.5/{1,2,4},
      134/0.5/{1,3,4}, 234/1.5/{2,3,4}}
      \node[subset] (s\id) at (\column,3) {$\elements$};

    \foreach \index in {1,2,3,4}
      \draw (s0) -- (s\index);
    \foreach \pair/\first/\second in {
      12/1/2, 13/1/3, 14/1/4, 23/2/3, 24/2/4, 34/3/4}
      \draw (s\first) -- (s\pair) (s\second) -- (s\pair);
    \foreach \triple/\first/\second/\third in {
      123/12/13/23, 124/12/14/24, 134/13/14/34, 234/23/24/34}
      \draw (s\first) -- (s\triple)
            (s\second) -- (s\triple)
            (s\third) -- (s\triple);
  \end{tikzpicture}
  \caption{Cover graph of the proper subsets of $\{1,2,3,4\}$.}
  \label{fig:proper-subset-cover-graph}
\end{figure}

Unlike simple graphs, graphs with loops can share FVC-TAR graphs with each other and even with simple graphs. The graphs in Figure \ref{fig:three-small-graphs} all have isomorphic FVC-TAR graphs, namely the graph in Figure \ref{fig:proper-subset-cover-graph}.

% Far-right graph in the top row.

\begin{proposition}\label{parallelclasses}
    Let $G$ be a graph, and let $G'$ be obtained from $G$ by replacing each parallel-edge class with a single edge. Then $\tau(G)=\tau(G')$ and $G$ and $G'$ have exactly the same vertex covers and failed vertex covers.
\end{proposition}
\begin{proof}
    A vertex set covers one edge of a parallel-edge class if and only if it covers every edge of that class. This includes parallel loops.
\end{proof}

\begin{proposition} \label{Loopy Stuff}
    If $L$ is a set of vertices with loops, then $C$ is a vertex cover of $G$ if and only if $L\subseteq C$ and $C\setminus L$ covers $G\setminus L$. Hence $\tau(G)=|L|+\tau(G\setminus L)$.
\end{proposition}
\begin{proof}
    For the forward direction, suppose that $C\subseteq V(G)$ is a vertex cover for $G$. Then, because loops are edges with only one incident vertex, $L \subseteq C$. Every edge of $G\setminus L$ has both endpoints outside of $L$, so its endpoint in $C$ must belong to $C\setminus L$. Thus $C\setminus L$ covers $G\setminus L$. 

    Conversely, choose $C$ so that $L \subseteq C$ and $C \setminus L$ covers $G \setminus L$. We partition $E(G)$ into two disjoint sets, those edges in $G \setminus L$ and those edges incident to a vertex in $L$. Because $L \subseteq C$, we know that $C$ will cover all the edges incident to some vertex in $L$. Additionally, because $C\setminus L$ covers $G \setminus L$, we know that $C$ will also cover those edges not incident to any vertex in $L$. Thus, $C$ is a vertex cover of $G$.

    Consequently, the vertex covers of $G$ are exactly the sets $L\cup D$ where $D$ is a vertex cover of $G\setminus L$. Taking minimum cardinalities gives $\tau(G)=|L|+\tau(G\setminus L)$.
\end{proof}

\begin{corollary}\label{cor-parallel}
    Let $G$ be a graph without parallel edges, and let $G'$ be obtained from $G$ by deleting each nonloop edge incident to a looped vertex, while retaining all vertices. Then $\tau(G)=\tau(G')$ and $G$ and $G'$ have exactly the same vertex covers and failed vertex covers.
\end{corollary}

\begin{proof}
    Let $L$ be the set of looped vertices. The two graphs have the same looped vertices and the same induced graph on $V(G)\setminus L$. Proposition \ref{Loopy Stuff} therefore gives exactly the same condition for a set to cover either graph.
\end{proof}

\begin{definition}Let \define{$R(G)$} be the graph with vertex set $V(G)$ obtained from $G$ by replacing each parallel-edge class by a single edge, and then removing each nonloop edge incident to a looped vertex. Thus one loop remains at every looped vertex. Note that $R(G)$ may have isolated vertices. If $G\cong R(G)$, then we call $G$ \define{reduced}.
\end{definition}

\begin{proposition}\label{reducedcong}
    Let $G$ be a graph. Then $\fvctar{G}\cong\fvctar{R(G)}$.
\end{proposition}

\begin{proof}
    The first step in the reduction preserves the vertex cover sets by Proposition \ref{parallelclasses} and the second step preserves them by Corollary \ref{cor-parallel}. Thus $G$ and $R(G)$ have the same failed vertex covers on the same vertex set. Adjacency follows using an argument similar to Lemma \ref{vertex count}, that is, by toggling one vertex at a time.
\end{proof}

We now extend Lemma \ref{findcenter} to the graphs considered in this section. In that lemma, which applies only to simple graphs with no isolated vertices, we showed that there is a unique center in the FVC-TAR graph. In the present section, that picture changes and there can be multiple centers. In particular, if there are no loops, the centers are exactly all subsets of the isolated vertices. If we allow loops, the centers are exactly the independent sets of the subgraph induced by the nonloop vertices; these sets may include any subset of the isolated vertices. Fortunately, the following proposition identifies a nonempty family of centers, namely subsets of the isolated vertices, each of which can be taken to the empty-set vertex by an automorphism of the FVC-TAR graph. In this way, we can still find ``the'' empty-set vertex of an unlabeled FVC-TAR graph.

\begin{proposition} \label{non-simple center vertex}
    Let $G$ be a graph with at least one edge, allowing loops, parallel edges, and isolated vertices. Then the empty set is a center of $\fvctar{G}$ and 
    \[ |V(G)|=\begin{cases}
        r(\fvctar{G})+2\quad&\textrm{ if }G\textrm{ has no loops,}\\
        r(\fvctar{G})+1\quad&\textrm{ if }G\textrm{ has a loop.}
    \end{cases}\] Let $C$ be the center of $\fvctar{G}$. Every vertex of maximum degree in the subgraph of $\fvctar{G}$ induced by $C$ can be taken to the empty-set vertex by an automorphism of $\fvctar{G}$. Further, these centers of maximum degree are exactly the subsets of the isolated vertices of $R(G)$.
\end{proposition}

\begin{proof}
    Observe first that, by Proposition \ref{reducedcong}, we may assume that $G$ is reduced. Also, the distance formula in Lemma \ref{symmdistance} still holds with loops and parallel edges allowed.

    Suppose $G$ has no loops. Let $I$ be the set of isolated vertices. Since $G$ is reduced and loopless, it is simple. By hypothesis, it also has an edge, and thus it follows that $G\setminus I$ is also simple, reduced, and has an edge. In addition, $G\setminus I$ has no isolated vertices and thus Lemma \ref{findcenter} applies. Every failed vertex cover of $G$ has the form $A\sqcup B$ where $A$ is a failed vertex cover of $G\setminus I$ and $B\subseteq I$. Let $A\cup B$ and $A'\cup B'$ be two such failed vertex covers. Then, Lemma \ref{symmdistance} states that \[d_{\fvctar{G}}(A\cup B,A'\cup B') = |(A\cup B)\ominus (A'\cup B')|=|A\ominus A'|+|B\ominus B'|=d_{\fvctar{G}}(A,A')+d_{\fvctar{G}}(B,B').\] For any fixed $B\subseteq I$, the quantity $|B\ominus B'|$ is maximized by choosing $B'=I\setminus B$, giving the value $|I|$. Thus $\epsilon_{\fvctar{G}}(A\cup B)=\epsilon_{\fvctar{G\setminus I}}(A)+|I|$.

    By Lemma \ref{findcenter}, $\fvctar{G\setminus I}$ has the unique center $\varnothing$ and radius $|V(G\setminus I)|-2$. It follows that $\fvctar{G}$ has radius $|V(G\setminus I)|-2+|I|=|V(G)|-2$. The sets $A\cup B$ that achieve this radius are exactly those sets satisfying $\epsilon_{\fvctar{G\setminus I}}(A)=|V(G\setminus I)|-2$, which by Lemma \ref{findcenter} is exactly when $A=\varnothing$. Since there are no restrictions on $B\subseteq I$, the centers of $\fvctar{G}$ are precisely the subsets of $I$. Within the subgraph induced by these centers, each center has exactly $|I|$ neighbors, obtained by adding or removing one element of $I$. Then each has maximum degree in that subgraph. Given one such center $Z\subseteq I$, define $f(A\cup B) = (A\cup B)\ominus Z=A\cup(B\ominus Z)$ for all $A\cup B\in V(\fvctar{G})$. Toggling elements of $I$ maintains the failed vertex cover property and adjacency, and $f$ is its own inverse. Hence $f$ is an automorphism of $\fvctar{G}$ with $f(Z)=\varnothing$.

    Now, suppose that $G$ has a loop and let $L$ be the set of looped vertices. A failed vertex cover has at most $|V(G)|-1$ vertices, and omitting a looped vertex attains this bound. Thus $\epsilon_{\fvctar{G}}(\varnothing)=|V(G)|-1$. We will establish that $|V(G)|-1$ is the minimum eccentricity among the vertices of $\fvctar{G}$. Let $S$ be a failed vertex cover. If $S$ contains an edge, then its complement is a failed vertex cover and the distance formula in Lemma \ref{symmdistance} gives $\epsilon_{\fvctar{G}}(S)=|V(G)|$. If $S$ does not contain an edge, then it is an independent set and its complement is a vertex cover. Since $S\cap L=\varnothing$, choose $l\in L$ and we have that $V(G)\setminus (S\cup \{l\})$ is a failed vertex cover since it leaves the loop at $l$ uncovered. Then Lemma \ref{symmdistance} provides that $d_{\fvctar{G}}(S,V(G)\setminus (S\cup \{l\}))=|S\ominus (V(G)\setminus (S\cup \{l\}))|=|V(G)|-1$. Thus $\epsilon_{\fvctar{G}}(S)=|V(G)|-1$ and it follows that $r(\fvctar{G})=|V(G)|-1$. Further, the centers of $\fvctar{G}$ are precisely the independent sets of $G\setminus L$.

    Let $C$ be the set of all such centers and let $K$ be the subgraph of $\fvctar{G}$ induced by $C$. For $S\in C$, let $N_{G\setminus L}(S)$ be the union of the neighborhoods in $G\setminus L$ of its elements. The $K$-neighbors of $S$ can be obtained from $S$ by removing an element of $S$ or by adding any vertex outside of $S\cup N_{G\setminus L}(S)$. Thus \[\deg_K(S)=|S|+|V(G\setminus L)\setminus(S\cup N_{G\setminus L}(S))|=|V(G)|-|L|-|N_{G\setminus L}(S)|.\] This degree is maximum exactly when $S$ consists entirely of isolated vertices of $G\setminus L$, that is, when $|N_{G\setminus L}(S)|=0$. Using the same automorphism definition and argument as the loopless case, we have the conclusion.
\end{proof}

Deciding what information can be recovered from an FVC-TAR graph is the center of what is shown in the following lemma. We require graphs have at least one edge because all edgeless graphs have the same - namely the empty - FVC-TAR graph, regardless of order. To see the reason for distinguishing between graphs with loops and those without, consider a graph consisting of one looped vertex and a graph consisting of two vertices with an edge between them. These graphs have isomorphic FVC-TAR graphs - namely $K_1$ - though their reduced graphs are not isomorphic. The three graphs in Figure \ref{fig:three-small-graphs} give another example. Each has the FVC-TAR graph shown in Figure \ref{fig:proper-subset-cover-graph}.

\begin{lemma} \label{The Big One}
    Let $G$ and $G'$ be graphs where either both have at least one loop or both have no loops and at least one edge. Then $\fvctar{G} \cong \fvctar{G'}$ if and only if $R(G) \cong R(G')$.
\end{lemma}
\begin{proof}
    By Proposition \ref{reducedcong}, $\fvctar{G}\cong\fvctar{R(G)}$ and $\fvctar{G'}\cong\fvctar{R(G')}$. Thus the reverse implication follows immediately. We then have only to prove that if $\fvctar{R(G)}\cong\fvctar{R(G')}$ then $R(G)\cong R(G')$.  Let $\phi:\fvctar{R(G)}\to\fvctar{R(G')}$ be an isomorphism. If both reduced graphs have no loops and no isolated vertices, then Proposition \ref{disconnected main theorem} already gives $R(G)\cong R(G')$. We now allow loops and isolated vertices.

    If $G$ and $G'$ both have at least one loop, then the same is true for $R(G)$ and $R(G')$. Proposition \ref{non-simple center vertex} then says that $|V(R(G))|=r(\fvctar{R(G)})+1=r(\fvctar{R(G')})+1=|V(R(G'))|$ and also that we can find a vertex in $\fvctar{R(G)}$ that can be treated as $\varnothing$ up to automorphism. We therefore assume that the isomorphism $\phi:\fvctar{R(G)}\to\fvctar{R(G')}$ satisfies $\phi(\varnothing)=\varnothing$. We can then use an argument similar to that in the proof of Proposition \ref{disconnected main theorem}. 
    
    Let $U\subseteq V(R(G))$ and $U'\subseteq V(R(G'))$ be the vertices whose singletons are failed vertex covers. Since these are exactly the neighbors of $\varnothing$, they can also be identified in the FVC-TAR graphs and $|U|=|U'|$ due to the isomorphism. We thus define a bijection $f:U\to U'$ by $f(v)=x$ when $\phi(\{v\})=\{x\}$. For $S\subseteq U$, define $f(S)=\{f(v):v\in S\}$. The rest of the proof follows exactly as in Proposition \ref{disconnected main theorem}, and is thus omitted. If $G$ and $G'$ are both loopless, then the same argument applies except that Proposition \ref{non-simple center vertex} says that $|V(R(G))|=r(\fvctar{R(G)})+2=r(\fvctar{R(G')})+2=|V(R(G'))|$.
    \end{proof}
    
All that remains to fully characterize the graphs with isomorphic FVC-TAR graphs is to consider graphs where one is loopless and the other is not. The following proposition covers this case. For $l\geq 1$, let $L_l$ denote the graph on $l$ vertices having one loop at each vertex and no other edges. Given two graphs $G$ and $H$, the \define{direct sum} $G\oplus H$ is the graph whose vertex set is $V(G)\cup V(H)$ and whose edge set is the disjoint union $E(G)\sqcup E(H)$. The graph $iK_1$ consists of $i$ isolated vertices.

First, we provide a lemma from \cite{ZFpaper} that was adapted from \cite[Theorems 2.23 and 8.13]{XTAR-YTAR} and included here for completeness. Note that failed vertex covering is a sub $Y$-set parameter and thus this lemma applies.

\begin{lemma}
\label{MaximalIFFdeg}
    Let $Y(G)$ be a sub $Y$-set parameter, $G$ be a graph of order $n$, and $T$ be a $Y$-set. Then, $\deg_{\mathfrak{Y}(G)}(T) \geq |T|$. Furthermore, $T$ is a maximal $Y$-set if and only if $\deg_{\mathfrak{Y}(G)}(T)=|T|$. 
\end{lemma}

\begin{proposition}\label{blahblahmixed}
    Let $G$ be loopless with at least one edge, and let $G'$ have at least one loop. Then $\fvctar{G}\cong\fvctar{G'}$ if and only if there exist integers $l\geq 1$ and $i\geq 0$ such that $R(G)\cong S_l\oplus iK_1$ and $R(G')\cong L_l\oplus iK_1$.
\end{proposition}

\begin{proof}
    Using Proposition \ref{reducedcong}, we assume both graphs are reduced. Assume that $\phi:\fvctar{G}\to\fvctar{G'}$ is an isomorphism. From Proposition \ref{non-simple center vertex}, we have that $|V(G)|-2=|V(G')|-1$ and thus $|V(G)|=|V(G')|+1$. Using an argument similar to that of Lemma \ref{The Big One}, we assume $\phi(\varnothing)=\varnothing$. By Lemma \ref{symmdistance}, if $S$ is a failed vertex cover of $G$, then $|\phi(S)|=d_{\fvctar{G'}}(\varnothing,\phi(S))=d_{\fvctar{G}}(\varnothing,S)=|S|$. By Lemma \ref{MaximalIFFdeg}, a set $S$ is maximal if and only if $\deg_{\fvctar{G}}(S)=|S|$. Thus $\phi$ preserves both cardinality and maximality.

    Since $G$ is loopless, the failed vertex covers are exactly those sets that miss at least one edge in $G$. Thus a maximal failed vertex cover $S$ misses exactly one edge and has cardinality $|V(G)|-2$. Suppose $G'$ has a maximal failed vertex cover $T$ of cardinality $|V(G')|-2$. This is equivalent to supposing that $G'$ has a nonloop edge $uv$ so that $T=V(G')-\{u,v\}$. Then there is a maximal failed vertex cover $S$ of $G$ so that $\phi(S)=T$. However, by our earlier argument, this would mean $|V(G)|-2=|S|=|T|=|V(G')|-2$, contradicting that $|V(G)|=|V(G')|+1$. Thus $G'$ has no nonloop edges and therefore $G'\cong L_l\oplus iK_1$ for some $l\geq 1$ and $i\geq 0$ with $|V(G')|=l+i$.

    Let $j$ be the number of isolated vertices of $G$. From Proposition \ref{non-simple center vertex}, the central vertices of $\fvctar{G}$ are exactly the subsets of its isolated vertices. Thus there are $2^j$ central vertices. Similarly, the central vertices of $\fvctar{G'}$ are the $2^i$ subsets of its isolated vertices. Because the FVC-TAR graphs of $G$ and $G'$ are isomorphic, it follows that $i=j$.

    If $G'$ has exactly one loop, then $|V(G)|=|V(G')|+1=1+i+1=i+2$. Then $G$ without its isolated vertices is a simple graph on two vertices and at least one edge, which is necessarily $K_2=S_1$. Thus $G\cong S_1\oplus iK_1$, as desired. We now assume that $G'$ has at least two loops. Then every singleton of $G'$ is a failed vertex cover and so $\deg_{\fvctar{G'}}(\varnothing)=|V(G')|=|V(G)|-1$. Since $\phi(\varnothing)=\varnothing$, it follows that $\deg_{\fvctar{G}}(\varnothing)=|V(G)|-1$. Thus $G$ has a unique singleton vertex cover $\{c\}$. By Proposition \ref{starunique}, this can only occur if the non-isolated portion of $G$ is a star. Since $G$ has $i$ isolated vertices, the number of leaves on the star is $|V(G)|-i-1=|V(G')|+1-i-1=(l+i)+1-i-1=l$. Therefore $G\cong S_l\oplus iK_1$.

    For the converse, suppose that $R(G)\cong S_l\oplus iK_1$ and $R(G')\cong L_l\oplus iK_1$. A failed vertex cover of $R(G)$ omits the central vertex of the star, contains a proper subset of its leaves, and contains an arbitrary subset of the isolated vertices. A failed vertex cover of $G'$ contains a proper subset of the looped vertices and an arbitrary subset of the isolated vertices. To form an isomorphism, we identify leaves of the star with looped vertices and isolated vertices with isolated vertices. The conclusion follows from Proposition \ref{reducedcong}.
    
\end{proof}

We are now ready to prove the main result.

\maintheorem*

\begin{proof}
    A graph has no failed vertex covers if and only if it is edgeless. Thus an edgeless graph has an FVC-TAR graph isomorphic to that of another graph exactly when the other graph is also edgeless. Now suppose both graphs have at least one edge. If both are loopless or both have at least one loop, then Lemma \ref{The Big One} provides our conclusion. If exactly one graph has a loop, then Proposition \ref{blahblahmixed} provides it.
\end{proof}

\begin{aistatement}During the preparation of this work, the authors used ChatGPT for mathematical exploration and computational code. After using this tool, the authors reviewed and edited the content as needed and they take full responsibility for the content of the published article.
\end{aistatement}

\end{document}